\documentclass[12pt]{amsart}
\usepackage{amsmath}
\usepackage{amsthm}
\usepackage{amssymb}

\newtheorem{theorem}{Theorem}[section]
\newtheorem{lemma}[theorem]{Lemma}
\newtheorem{corollary}[theorem]{Corollary}

\theoremstyle{definition}
\newtheorem{definition}[theorem]{Definition}

\theoremstyle{remark}
\newtheorem{remark}[theorem]{Remark}

\numberwithin{equation}{section}

\newcommand{\be}{\begin{equation}}
\newcommand{\ee}{\end{equation}}
\newcommand{\F}{\mathbb{F}}
\newcommand{\Fp}{\mathbb{F}_p}
\newcommand{\Z}{\mathbb{Z}}

\begin{document}

\title{Primality tests for Fermat numbers and $2^{2k+1}\pm2^{k+1}+1$.}


\author{Yu Tsumura}
\address{Department of Mathematics, Purdue University
150 North University Street, West Lafayette, Indiana 47907-2067
}

\email{ytsumura@math.purdue.edu}
\thanks{}


\date{}

\dedicatory{}

\begin{abstract}
Robert Denomme and Gordan Savin made a primality test for Fermat numbers $2^{2^k}+1$ using elliptic curves.
We propose another primality test using elliptic curves for Fermat numbers  and also give primality tests for integers of the form $2^{2k+1}\pm2^{k+1}+1$.

\end{abstract}

\maketitle


\section{Introduction.}

The integers of the form $2^{2^k}+1$ with $k\geq 0$ are called Fermat numbers, named after  Pierre de Fermat.
For $k=0$, $1$, $2$, $3$, $4$, Fermat numbers are prime.
Fermat conjectured that all numbers of this form were prime numbers. 
However, in 1732 Leonhard Euler disproved this conjecture by factoring the fifth Fermat number $2^{2^5}+1=641\cdot 6700417$. 
Not only  was it disproved, but also  no other Fermat primes  have been discovered when $k>4$.
So checking the primality or finding factors of Fermat numbers attracts many people.

Let us define the notation used in this paper.
\begin{definition}
Let $F_k=2^{2^k}+1$, $G_k=2^{2k+1}+2^{k+1}+1$, and $H_k=2^{2k+1}-2^{k+1}+1$, where $k$ is assumed to be a positive integer.
$F_k$ is called the $k$th Fermat number.
\end{definition}

In 1877, Pepin gave a very efficient  primality test for Fermat numbers.
\begin{theorem}\label{thm:pepin}{(Pepin test).}
For $k\geq1$, $F_k=2^{2^k}+1$ is prime if and only if $3^{(F_k-1)/2}\equiv -1 \pmod {F_k}$.
\end{theorem}
\begin{proof}
See Theorem 4.1.2 in \cite{Crandall-Pomerance}.
\end{proof}

In this paper, we study group structures of  elliptic curves defined over finite fields of order $F_k$, $G_k$, and $H_k$ (if they are prime). The essential role is the action of an endomorphism $[1+i]$ on the curves.
After that we  use the information of the group structure to give two primality tests for Fermat numbers which can be regarded as an elliptic version of the Pepin test.
Also, we give similar results for integers of the form $2^{2k+1}\pm 2^{k+1}+1$.

The original work in this direction was done by Benedict H. Gross in \cite{Gross} for Mersenne numbers and by Robert Denomme and Gordan Savin in \cite{Denomme-Savin} for Fermat numbers and integers of the form $3^{2^k}-3^{2^{k-1}}+1$ and $2^{2^k}-2^{2^{k-1}}+1$, where $k$ is a positive integer.
Gross used the formula of the multiplication by 2 as a recursive formula and Denomme and Savin used the formula of the action of $[1+i]$ as a recursive formula for Fermat numbers.
In this paper, we obtain the same primality test as Denomme and Savin in a slightly different approach and also give a new primality test which uses the formula of the multiplication by 2 for Fermat numbers. Also, by the same method we give new primality tests for $G_k$, $H_k$.
As you notice by the following proofs, $F_k$, $G_k$ and $H_k$ are the only numbers to which this method applies.

We saw in Theorem \ref{thm:pepin} that there is a fast primality test for $p=F_k$.
There are also fast primality tests for $p=G_k$ and $p=H_k$.
For example, one could use Corollary 1 or Theorem 5 of \cite{Brillhart-Lehmer-Selfridge}.
These tests apply because $p-1$ is divisible by a power of $2$ near $\sqrt{p}$.
These tests determine the primality of $p$ of these three special forms in polynomial time.
Our new tests below also run in polynomial time and are the first such tests using elliptic curves.

\section{Group Structure.}

The next theorem allows us to determine the order of certain elliptic curve groups.

\begin{theorem}\label{thm:group-structure}
Let $p\equiv 1 \pmod 4$ be an odd prime and let $m \not\equiv 0 \pmod p$ be a fourth power mod $p$.
Let $E$ be an elliptic curve defined by $y^2=x^3-mx$.
Let $p=a^2+b^2$, where $a$, $b$ are integers with $b$ even and $a+b\equiv1 \pmod 4$.
Let $E(p)$ be the elliptic curve $E$ defined over $\Fp$.
Then we have $\#E(p)=p+1-2a$.
\end{theorem}
\begin{proof}
See Theorem 4.23, page 115 in \cite{Washington}.
\end{proof}

From now on, we fix an elliptic curve $E:y^2=x^3-mx$, where $m\not\equiv 0 \pmod p$ is a fourth power mod a prime $p$.
We  denote by $E(p)$  the elliptic curve group $E$ defined over finite field $\Fp$ when $p$ is prime.
Also let $E(\bar{\F}_p)$ be the elliptic curve $E$ defined over the algebraic closure $\bar{\F}_p$ of $\Fp$ and we  denote by $E[n]$ the elements in $E(\bar{\F}_p)$ whose orders divide $n$.
\begin{corollary}\label{cor:order}
\begin{enumerate}
	\item If $F_k$ is prime, then $\#E(F_k)=2^{2^k}$.
	\item If $G_k$ is prime, then $\#E(G_k)=2^{2k+1}$.
	\item If $F_k$ is prime, then $\#E(H_k)=2^{2k+1}$.
\end{enumerate}
\end{corollary}
\begin{proof}
Let us first consider $F_k$.
The decomposition into two squares is 
$F_k=2^{2^k}+1=1^2+(2^{2^{k-1}})^2$ and $1+(2^{2^{k-1}})\equiv 1 \pmod 4$.
Hence by Theorem \ref{thm:group-structure}, $\#E(F_k)=F_k+1-2=2^{2^k}$.

Next, let $a=2^k+1$ and $b=2^k$.
Then we have $G_k=a^2+b^2$ and $a+b\equiv 1 \pmod 4$.
Hence we have $\#E(G_k)=G_k+1-2(2^k+1)=2^{2k+1}$ by Theorem \ref{thm:group-structure}.

Similarly, let $a=-(2^k-1)$ and $b=2^k$.
Then we have $H_k=a^2+b^2$ and $a+b\equiv 1 \pmod 4$.
Hence $\#E(H_k)=H_k+1+2(2^k-1)=2^{2k+1}$.

\end{proof}

The next lemma gives  information on the group structures of $E(p)$ and $E[n]$.
\begin{lemma}\label{lem:basic}
Let $E$ be an elliptic curve over a finite field $\Fp$. 
Then we have
\[E(p) \cong \Z_{n_1} \oplus \Z_{n_2}\]
for some positive integers $n_1$ and $n_2$ with $n_1 |n_2$.
Also, if $n$ is a positive integer which is not divisible by $p$,  then we have
\[E[n] \cong \Z_n \oplus \Z_n.\]
\end{lemma}
\begin{proof}
See Theorem 3.1 and Theorem 4.1 in \cite{Washington}.
\end{proof}

Let $p$ denote one of $F_k$, $G_k$ and $H_k$.
Suppose $p$ is  prime.
By Corollary \ref{cor:order} and Lemma \ref{lem:basic}, the group structure is  $E(p)\cong \Z_{2^\alpha} \oplus \Z_{2^\beta}$ with $\alpha \leq \beta$ and $\alpha+\beta=2^k$ if $p=F_k$ and $\alpha+\beta=2k+1$ if $p=G_k$ or $p=H_k$.
Since $m$ is a 4th power, all the roots of $x^3-mx$ are in $\Fp$ and also in the subgroup $E[2]\cong \Z_2 \oplus \Z_2$ by Lemma \ref{lem:basic}.
Then $\Z_2 \oplus \Z_2\cong E[2]\subset E(p)$, hence $E(p)$ is not cyclic. 
However, we can determine the group structure of $E(p)$ precisely.
First we need two lemmas.
\begin{lemma}\label{lem:subset}
Let $n$ be a positive integer which is not divisible by a prime $p$.
Let $\phi$ be the Frobenius endomorphism on $E(\bar{\F}_p)$ given by $\phi(x,y)=(x^p,y^p)$.
Then $E[n]\subset E(p)$ if and only if  $\phi -1$ is divisible by $n$ in ${\rm End}(E)$.
\end{lemma}
\begin{proof}
See Lemma 1 in \cite{Ruck}.
\end{proof}

\begin{lemma}\label{lem:frob}
If $\#E(p)=p+1-A$, then the Frobenius endomorphism $\phi$ satisfies $\phi^2-A\phi+p=0$ as an endomorphism of $E$.
\end{lemma}
\begin{proof}
See Theorem 4.10, page 101 in \cite{Washington}.
\end{proof}

\begin{theorem}\label{thm:fermatgp}
Suppose $F_k$ is prime.
Then we have
\[E(\F_k)\cong \Z_{2^{2^{k-1}}} \oplus \Z_{2^{2^{k-1}}}. \]
\end{theorem}
\begin{proof}

Since $\#E(F_k)=F_k+1-2$, the Frobenius endomorphism $\phi$ satisfies $\phi^2-2\phi+F_k=0$ in ${\rm End}(E)$ by Lemma \ref{lem:frob}, and hence $(\phi-1)^2=-2^{2^k}$. 
Since ${\rm End}(E)\cong \Z[i]$ (see chapter 10 in \cite{Washington} ), it is a unique factorization domain. 
Therefore $\phi-1=\pm i2^{2^{k-1}}$, and hence $2^{2^{k-1}}$ divides $\phi-1$.
Then $E[2^{2^{k-1}}]\subset E(F_k)$ by Lemma \ref{lem:subset}.
Since $E[2^{2^{k-1}}]\cong \Z_{2^{2^{k-1}}} \oplus \Z_{2^{2^{k-1}}}$ by Lemma \ref{lem:basic}, we have $\#E[2^{2^{k-1}}]=(2^{2^{k-1}})^2=2^{2^k}=\#E(F_k)$.
Therefore we have $E(F_k)=E[2^{2^{k-1}}]\cong \Z_{2^{2^{k-1}}} \oplus \Z_{2^{2^{k-1}}}$.
\end{proof}
\begin{theorem}\label{thm:Gkgp}
Suppose $G_k$ is prime. 
Then we have 
\[E(G_k)\cong \Z_{2^k}\oplus\Z_{2^{k+1}}.\]
\end{theorem}
\begin{proof}
From Corollary \ref{cor:order}, we know that $\#E(G_k)=2^{2k+1}=G_k+1-2(2^k+1)$.
Hence the Frobenius endomorphism $\phi$ satisfies $\phi^2-2(2^k+1)\phi+G_k=0$.
Then we have $(\phi-1)^2-2^{k+1}(\phi-1)+2^{2k+1}=0$.
Therefore, $\phi-1=2^k(1\pm i)$.
Hence $2^k$ divides $\phi-1$ and we have $E[2^k]\subset E(G_k)$ by Lemma \ref{lem:subset}.
Since $\#E[2^k]=2^{2k}$ and $\# E(G_k)=2^{2k+1}$, the group structure of $E(G_k)$ must be  $E(G_k)\cong  \Z_{2^k}\oplus\Z_{2^{k+1}}$ by Lemma \ref{lem:basic}.
\end{proof}
\begin{theorem}\label{thm:Hkgp}
Suppose $H_k$ is prime. 
Then we have 
\[E(H_k)\cong \Z_{2^k}\oplus\Z_{2^{k+1}}.\]
\end{theorem}
\begin{proof}
Just note that the Frobenius endomorphism satisfies $\phi^2+2(2^k-1)\phi+H_k=0$.
Hence $\phi-1=(-1\pm i)2^k$. The rest of the proof is identical to that of Theorem \ref{thm:Gkgp}.
\end{proof}

\section{Primality test}
Again let $p$ be one of $F_k$, $G_k$ and $H_k$.
As we noted in the proof of Theorem \ref{thm:fermatgp}, $E$ has  complex multiplication by $\Z[i]$.
For a detailed explanation about complex multiplication, see chapter 10 in \cite{Washington}.
The action of $i$ on $(x,y)\in E$ is given by $[i] \cdot (x,y)=(-x,iy)$, where the $i$ in $(-x,iy)$ is a 4th root of unity in $\Fp$. 
This $i$ exists in $\Fp$ since $p\equiv 1 \pmod 4$.
Note that as an endomorphism, $i$ has degree $1$ and hence it is an isomorphism.
Now, let us denote $\eta=1+i$ in ${\rm End}(E)$.
This endomorphism is very important in this paper.
Let us describe the action of $\eta$ on $(x,y)$ explicitly.
Let $\eta \cdot(x,y)=(x',y')$. 
We have 
\begin{align*}
\eta\cdot(x,y)&=[1+i]\cdot(x,y)=(x,y)+[i]\cdot(x,y)=(x,y)+(-x,iy)
\end{align*}
and by the elliptic curve addition, it is equal to
\begin{equation}\label{eq:eta1}
\left( \left(\frac{(1-i)y}{2x}\right)^2, y' \right)
\end{equation}
\begin{equation}\label{eq:eta2}
=\left(\frac{x^2-m}{2ix} ,y' \right),
\end{equation}
where $y'=\left(\frac{(1-i)y}{2x} \right) (x-x')-y$.
Note that by the  equation (\ref{eq:eta1}), the $x$-coordinate $x'$ of  $\eta \cdot(x,y)$ is a square and by the equation (\ref{eq:eta2}), $x'$ can be computed without using $y$.
Also note that $\eta$ has degree $2$, hence $\# {\rm Ker}(\eta)=2$. 
Clearly, $(0,0)$ is in the kernel and so ${\rm Ker}(\eta)=\{\infty, (0,0) \}$, where $\infty$ is the identity of $E$.

Note that $\eta^2=2i$ and $\eta^{2l}=\epsilon2^l$, where $l$ is a positive integer and $\epsilon=\pm 1$, $\pm i$.
Since $\epsilon=\pm 1$, $\pm i$ are isomorphism, we do not care about this factors.
We will use $\epsilon$ for $\pm 1$, $\pm i$ in this paper, but $\epsilon$ might have different values at each occurrence.

\subsection{Primality test for Fermat numbers.}
Now we can state a theorem which can be converted into a primality test.
\begin{theorem}\label{thm:fermatmain}
Let $\eta=1+i$ in ${\rm End}(E)$. Let $P=(x,y)$ on $E$, where $x$ is a quadratic non-residue mod $F_k$.
Then $F_k$ is prime if and only if $\eta^{2^{k}-1}P=(0,0)$.
\end{theorem}
\begin{proof}
Suppose $F_k$ is prime.
In the proof of  Theorem \ref{thm:fermatgp}, we have seen that $\phi-1=\epsilon 2^{2^{k-1}}=\epsilon \eta^{2^k}$.
Hence, we have ${\rm Ker}(\eta^{2^k})={\rm Ker}(\phi-1)=E(F_k)$.
Since $\#{\rm Ker}(\eta)=2$ and $\#E(F_k)=2^{2^k}$, we have ${\rm Ker}(\eta^s)={\rm Im}(\eta^{2^k-s})$ for $s=1$, $2$, \ldots, $2^k$.
Assume $P=\eta Q$ for some $Q\in E(F_k)$.
Then as we noted above, the $x$-coordinate $x$ of $\eta Q=P$ is a square.
However, we assumed that $x$ is a quadratic non-residue mod $F_k$, hence $P$ is not in the image of $\eta$.
Observe that $\eta^{2^k-1} P \neq \infty$ since otherwise $P\in {\rm Ker}(\eta^{2^k-1})={\rm Im}(\eta)$, but $P \notin {\rm Im}(\eta)$.
Since  $\eta^{2^k-1} P \neq \infty$ and $\eta^{2^k}=\infty$, we have $\eta^{2^k-1} P=(0,0)$.

Conversely, suppose $\eta^{2^k-1}P=(0,0)$.
Assume $F_k$ is composite and let $q$ be a prime divisor such that $q\leq \sqrt{F_k}$.
It is known that a divisor of a Fermat number is congruent to $1$ modulo $4$.
(See \cite{Crandall-Pomerance}).
Then $\eta^{2^k-1}P=(0,0)$ holds in the reduction $E(q)$.
It follows that $2^{2^{k-1}-1}P=\epsilon \eta^{2^k-2}P\neq \infty$.
Also we have $2^{2^{k-1}}P=\epsilon\eta^{2^{k}}P=\infty$, therefore $P$ has order $2^{2^{k-1}}$.
Assume that $\{ P, iP \}$ is a basis of $E[2^{2^{k-1}}]$.
Note that $iP\in E(q)$ since   $i\in \F_{q}$ when $q\equiv1 \pmod4$.
So we have $E[2^{2^{k-1}}] \subset E(q)$, hence $2^{2^k} \leq \#E(q)$.
However, $\#E(q)\leq (\sqrt{q}+1)^2$ by Hasse's Theorem.
Hence, we have $q^2-1\leq F_k^2-1=2^{2^k}\leq \#E(q)\leq (\sqrt{q}+1)^2 $.
This inequality holds only for $q=2$.
However, clearly $q$ is an odd prime.
Hence it is a contradiction.
Therefore $F_k$ is prime.

To complete the proof, we need to prove that $\{ P, iP \}$ is a basis of $E[2^{2^{k-1}}]$.
Suppose $uP+v(iP)=\infty$ for some integers $u$, $v$.
Let $u=2^{\alpha}u'$ and let $v=2^{\beta}v'$ with $u'$, $v'$ odd.
Since the order of $P$ is a power of $2$, we have $\alpha=\beta$.
Now $(u'+v'i)(2^{\alpha}P)=\infty \Rightarrow (u'^2+v'^2)(2^{\alpha}P)=\infty \Rightarrow u'^2+v'^2\equiv 0 \pmod {2^{2^{k-1}-\alpha}}$.
Since $u'^2+v'^2\equiv 2 \pmod 4$, the above congruence holds only if $\alpha=2^{k-1}$ or  $\alpha=2^{k-1}-1$.
If $\alpha=2^{k-1}$, then $u\equiv v \equiv 0 \pmod {2^{2^{k-1}}}$, and hence they are independent.

Next let us consider the case $\alpha=2^{k-1}-1$.
Let $P'=(2^{2^{k-1}-1})P$. Then $P'$ has order $2$.
Hence $P'$ is either $(0,0)$ or $(\pm \sqrt{m}, 0)$.
However, $\eta P'=\eta\cdot (\epsilon \eta^{2^k-2})P=\epsilon\eta^{2^k-1}P\neq \infty$, hence we have $P'\neq (0,0)$.
Therefore, $P'$ is either $(\sqrt{m}, 0)$ or $(-\sqrt{m}, 0)$.
If $P'=(\sqrt{m}, 0)$, then $\infty=(u'+v'i)(\sqrt{m}, 0)=u'(\sqrt{m}, 0)+v'(-\sqrt{m}, 0)$ with odd $u'$, $v'$.
Since $\{(\sqrt{m}, 0), (-\sqrt{m}, 0) \}$ is a basis for $E[2]$, they cannot be dependent with odd coefficients.
The same thing happens when $P'=(-\sqrt{m}, 0)$.
Therefore,  $P$ and $iP$ are independent, and this completes the proof.
\end{proof}

Hence, to check the primality of Fermat numbers, we need to calculate $\eta^{2^{k}-1}P$ for a point $P$ with a quadratic non-residue $x$-coordinate mod $F_k$.
However, we need not to calculate a $y$-coordinate since when an $x$-coordinate is $0$, so is the $y$-coordinate.
Also as noted above, to calculate the $x$-coordinate of $\eta P$, the $y$-coordinate of $P$ is not used.

For example, take $m=1$ and $P=(5,2\sqrt{30})$ on $E:y^2=x^3-x$.
It is straightforward to check $5$ is a quadratic non-residue and $30$ is a quadratic residue mod $F_k$.
Hence $P$ satisfies the conditions of Theorem \ref{thm:fermatmain}.

Here is the algorithm to check the primality for $F_k$.
Let $x_0=5$ and let 
\[x_j= \frac{x_{j-1}^2-1}{2ix_{j-1}} \]
if $\gcd(x_{j-1}, F_k)=1$ for $j\geq 1$.
Note that $x_j$ is the $x$-coordinate of $\eta^jP$.
Here $i$ is a primitive 4th root of unity in $F_k$ and it is explicitly $i=2^{2^{k-1}}$.
If $\gcd(x_j, F_k)>1$ for some $j<2^k-1$, then $F_k$ is composite and we terminate the algorithm.
If we calculate $x_{2^k-1}$ and it is $0$, then $F_k$ is prime.
If $x_{2^k-1}\neq 0$, then $F_k$ is composite.

\begin{remark}
We do not need to find $\sqrt{30}$ mod $F_k$ explicitly.
We just needed to know that the point $P=(5,2\sqrt{30})$ is on $E:y^2=x^3-x$.
What we need is  only the $x$-coordinate in the algorithm.
\end{remark}

An alternative primality test can be deduced by noting equivalent conditions as in the next lemma.
\begin{lemma}
Let $P$ be a point on $E$ with a quadratic non-residue $x$-coordinate mod $F_k$.
Then $\eta^{2^k-1}P=(0,0)$ if and only if $2^{2^{k-1}}P=(\sqrt{m},0)$ or $(-\sqrt{m}, 0)$.
\end{lemma}
\begin{proof}
Suppose $\eta^{2^k-1}P=(0,0)$.
Then we have $\eta(2^{2^{k-1}-1}P)=\epsilon\eta \cdot \eta^{2^k-2}P=(0,0) $.
Therefore we have $2^{2^{k-1}-1}P \neq \infty$, $(0,0)$, otherwise the image by $\eta$ is $\infty$.
Also, we have $2(2^{2^{k-2}}P)=2^{2^{k-1}}P=\epsilon \eta^{2^k}P=\epsilon \eta(0,0)=\infty$.
Therefore $2^{2^{k-1}}P \in E[2]\setminus \{\infty, (0,0) \}$.
That is, $2^{2^{k-1}}P=(\sqrt{m},0)$ or $(-\sqrt{m}, 0)$.

Conversely, suppose  $2^{2^{k-1}}P=(\pm\sqrt{m},0)$.
We have 
\[(0,0)=\eta(\pm\sqrt{m},0)=\eta(2^{2^{k-1}})P=\epsilon \eta^{2^k-1}P.\]
Hence, we have $\eta^{2^k-1}P=(0,0)$.
\end{proof}

So now we have shifted from the multiplication by $\eta$ to the multiplication by $2$.
Multiplication by $2$ of a point $P=(x,y)$ on the elliptic curve $E:y^2=x^3-mx$ is described as follow.
\[2(x,y)=\left(\frac{x^4+2mx^2+m^2}{4(x^3-mx)},yR(x) \right)\]
for some rational function $R(x)$.
(See Example 2.5, page 52 in \cite{Washington}.)
Let $P=(x_0,y_0)$ be a point on $E$ with a quadratic non-residue $x$-coordinate mod $p$.
Let 
\[x_j=\frac{x_{j-1}^4+2mx_{j-1}^2+m^2}{4(x_{j-1}^3-mx_{j-1})}\]
modulo $F_k$ if $\gcd((x_{j-1}^3-mx_{j-1}), F_k)=1$ for $j\geq 1$ inductively. 
Hence $x_j$ is the $x$-coordinate of $2^jP$.
If we can proceed to calculate $x_{2^{k-1}-1}$ and this is $\pm \sqrt{m}$, then $F_k$ is prime.
Otherwise $F_k$ is composite.

For example, let us consider the same example as above.
Let $m=1$ and $P=(5,2\sqrt{30})$ on E.
Then the algorithm to check the primality for $F_k$ is as follows.
Let $x_0=5$ and we define inductively 
\[x_j=\frac{x_{j-1}^4+2x_{j-1}^2+1}{4(x_{j-1}^3-x_{j-1})}\]
 if $\gcd((x_{j-1}^3-x_{j-1}), F_k)=1$ for $j\geq 1$.
 If $\gcd((x_{j-1}^3-x_{j-1}), F_k)=1$ for some $j<{2^{k-1}-1}$, then $F_k$ is composite and we terminate the algorithm.
If we  calculate $x_{2^{k-1}-1}$ and this is $\pm 1$, then $F_k$ is prime.
Otherwise $F_k$ is composite.

\begin{remark}
Although the recursion formula for $x_j$ looks more complicated than before, the number of recursions is reduced to $2^{k-1}-1$ from $2^k-1$.
\end{remark}
\subsection{Primality test for $2^{2k+1}+2^{k+1}+1$.}
\begin{theorem}\label{thm:Gk}
Let $P=(x,y)$ be a point on $E$, with $x$ is a quadratic non-residue mod $G_k$.
Then $G_k$ with $k\geq2$ is prime if and only if $\eta^{2k-1}P\in E[2]\setminus \{\infty\}$.
\end{theorem}
\begin{proof}
Suppose $G_k$ is prime.
We have  $\#(\eta^{2k}E(G_k))=\#(\epsilon2^kE(G_k))=2$.
We have seen that $\phi-1=\epsilon\eta2^k=\epsilon \eta^{2k+1}$ when $G_k$ is prime in the proof of Theorem \ref{thm:Gkgp}.
Since ${\rm Ker}(\phi-1)= E(G_k)$, we have $\eta(\eta^{2k}E(G_k))=\infty$, and therefore $\eta^{2k-1}E(G_k)=E[2]$.

Now that we know that $E(G_k)={\rm Ker}(\eta^{2k+1})$ and $\#{\rm Ker}(\eta)=2$ in addition to $\#E(G_k)=2^{2k+1}$, it is easy to see that ${\rm Ker}(\eta^s)={\rm Im}(\eta^{2k+1-s})$, for $ s=0$, $1$, \ldots, $2k+1$.
Since $x$ is not a square mod $p$, $P$ is not in the image of $\eta$.
Hence, we have $\eta^{2k-1}P\in E[2]\setminus \{\infty\}$.
Let us show this.
If $\eta^{2k-1}P=\infty$, then $P\in {\rm Ker}(\eta^{2k-1})={\rm Im}(\eta^2)$.
Since $P$ is not in the image of $\eta$, this is a contradiction.
Hence $\eta^{2k-1}P\neq \infty$.

Conversely, suppose $\eta^{2k-1}P\in E[2]\setminus \{\infty\}$.
Assume $G_k$ is composite and let $q$ be a prime divisor of $G_k$ such that $q\leq \sqrt{G_k}$.
Then $\eta^{2k-1}P\in E[2]\setminus \{\infty\}$ holds in the reduction $E(q)$.
Then $\eta^{2k-1}P$ is one of $(0,0)$ or $(\pm \sqrt{m},0)$.
If $\eta^{2k-1}P=(0,0)$, then we have $2^{k-1}P=\epsilon\eta^{2k-2}P\neq\infty$ and $2^kP=\epsilon\eta^{2k}P=\infty $.
Therefore $P$ has order $2^k$.
If $\eta^{2k-1}P=(\sqrt{m},0)$, then let $P'=\eta P$.
Then we have $\eta^{2k-1}P'=\eta(\sqrt{m}, 0)=(0,0)$. 
This is the same situation as the case $\eta^{2k-1}P=(0,0)$, hence $P'$ has order $2^k$. 
The case $\eta^{2k-1}P=(-\sqrt{m},0)$ is similar and $\eta P$ has order $2^k$.
We have seen in any case, there exists a point ($P$ or $\eta P$) of order $2^k$.
Let $R$ denote this point.
Let us assume that $\{R, iR\}$ is a basis for $E[2^k]$.
It is easy to check that every divisor of $G_k$ is congruent to $1$ modulo $4$.
So $iR \in E(q)$ and hence  $E[2^k]\subset  E(q)$.
Therefore we have
\[ 2^{2k}=\#E[2^k]\leq \# E(q)\leq (\sqrt{q}+1)^2\leq(G_k^{1/4}+1)^2. \]
However, this inequality does not hold for $k\geq 2$, and therefore $G_k$ is prime.

To complete the proof, we need to show that $\{R, iR\}$ is a basis for $E[2^k]$.
Suppose $uR+v(iR)=\infty$ for some integers $u$, $v$.
Let $u=2^{\alpha}u'$ and let $v=2^{\beta}v'$ with $u'$, $v'$ odd.
Since the order of $R$ is a power of $2$, we have $\alpha=\beta$.
Now $(u'+v'i)(2^{\alpha}R)=\infty \Rightarrow (u'^2+v'^2)(2^{\alpha}R)=\infty \Rightarrow u'^2+v'^2\equiv 0 \pmod {2^{k-\alpha}}$.
Since $u'^2+v'^2\equiv 2 \pmod 4$, the above congruence holds only if $\alpha=k$ or  $\alpha=k-1$.
If $\alpha=k$, then $u\equiv v \equiv 0 \pmod {2^k}$, and hence they are independent.

Next, let us consider the case $\alpha=k-1$.
Let $R'=2^{k-1}R$. 
Then $P'$ has order $2$.
Hence $R'$ is either $(0,0)$ or $(\pm \sqrt{m}, 0)$.
However, we have $\eta R'=\eta\cdot (\epsilon\eta^{2k-2})R =\epsilon\eta^{2k-1}R$
\begin{align*}
 = \left\{
\begin{array}{ll}
 \epsilon \eta^{2k-1}P\neq \infty & \textrm{ if } R=P \\
             \eta\cdot\eta^{2k-1}P=\eta(1\pm\sqrt{m}, 0)=(0,0) \neq \infty & \textrm{ if } R=\eta P.
\end{array}\right.
\end{align*}
Hence $R'\neq (0,0)$.
Therefore $P'$ is either $(\sqrt{m}, 0)$ or $(-\sqrt{m}, 0)$.
If $R'=(\sqrt{m}, 0)$, then $\infty=(u'+v'i)(\sqrt{m}, 0)=u'(\sqrt{m}, 0)+v'(-\sqrt{m}, 0)$ with odd $u'$, $v'$.
Since $\{(\sqrt{m}, 0), (-\sqrt{m}, 0) \}$ is a basis for $E[2]$, they cannot be dependent with odd coefficients.
The same thing happens when $R'=(\sqrt{m}, 0)$.
Therefore, $R$ and $iR$ are independent.
\end{proof}
To use Theorem \ref{thm:Gk}, we need to find a point on $E$ whose $x$-coordinate is a quadratic non-residue mod $G_k$.
It is straightforward to check the following.
\begin{itemize}
	\item $3$ is a quadratic non-residue mod $G_k$ if and only if $k$ is even. 
    \item $5$ is a quadratic non-residue mod $G_k$ if and only if $k\equiv 1 \pmod 4$.
  Also If $k\equiv0, 3 \pmod 4$, then $G_k$ is divisible by $5$.
  \item $7$ is a quadratic non-residue mod $G_k$ for all $k\geq 1$.
\end{itemize}
 
Using these facts, we can choose specific initial values depending on $k$.
Since  $G_k$ is composite when $k\equiv0, 3 \pmod 4$ from the above fact, we only need to consider the cases when $k\equiv 1 \pmod 4$ and $k\equiv 2 \pmod 4$.

When $k\equiv 2 \pmod 4$, we take $m=1$ and $P=(7,4\sqrt{21})$ on $E:y^2=x^3-x$.
Note that $21=3\cdot7$ is a quadratic residue mod $G_k$ since both $3$ and $7$ are quadratic non-residues.

When $k\equiv 1 \pmod 4$ and $k>1$, we can take $m=3^4$ ($3$ does not divide $G_k$) and $P=(5,2\sqrt{-70})$ on $E:y^2=x^3-3^4x$.
Note that $-70=-2\cdot 5 \cdot 7$ is a quadratic residue mod $G_k$ since $-2$ is a quadratic residue (because $G_k\equiv1 \pmod 8)$ and $5$ and $7$ are quadratic non-residues from the above facts.

Then the algorithm to check the primality of $G_k$ is as follows.
Let $x_0=7$ when $k\equiv 2 \pmod 4$ and $x_0=5$ when  $k\equiv 1 \pmod 4$.
Then let $x_j= (x_{j-1}^2-1)/(2ix_{j-1})$  if $\gcd(x_{j-1}, G_k)=1$ for $j\geq 1$ inductively.
As before this is the $x$-coordinate of $\eta^jP$.
If $\gcd(x_{j-1}, G_k)>1$ for some $j<2k-1$, then $G_k$ is composite and we terminate the algorithm.
If we  calculate $x_{2k-1}$ and this is $\pm 1$, then $G_k$ is prime.
Otherwise, $G_k$ is composite.

\subsection{Primality test for $2^{2k+1}-2^{k+1}+1$.}

Now let us discuss $H_k=2^{2k+1}-2^{k+1}+1$.
By Theorem \ref{thm:Hkgp}, we know that $\phi-1=\epsilon \eta^{2k+1}$.
Therefore the proof of the next theorem is identical to that of Theorem \ref{thm:Gk}.

\begin{theorem}\label{thm:Hk}
Let $P=(x,y)$ be a point on $E$, with $x$ is a quadratic non-residue mod $H_k$.
Then $H_k$, $k\geq2$ is prime if and only if $\eta^{2k-1}P\in E[2]\setminus \{\infty\}$.
\end{theorem}

Again to use Theorem \ref{thm:Hk}, we need to find a point on a curve whose $x$-coordinate is a quadratic non-residue mod $H_k$.
The following is easy to check.
\begin{itemize}
	\item $3$ is a quadratic non-residue mod $H_k$ if and only if $k$ is even.
	\item $5$ is a quadratic non-residue mod $H_k$ if and only if $k\equiv 3 \pmod 4$.
	Also when $k\equiv 1,2 \pmod 4$, $5$ divides $H_k$.
	\item When $k\equiv 4 \pmod {12}$, $13$ divides $H_k$.	
\end{itemize}

Hence when $k\equiv 3 \pmod 4$, we can take $m=1$ and a point $(5,2\sqrt{30})$ on $E: y^2=x^3-x$.
Here $30=2\cdot3\cdot 5$ is a quadratic residue by the above facts.

The remaining cases are when $k\equiv 0, 8 \pmod {12}$, otherwise $5$ or $13$ divides $H_k$.
However, it seems difficult to find a suitable  small initial value.
So we further divide the cases into $k\equiv 0$, $8$, $12$, $20$, $24$, $32$, $36$, $44$ $\pmod {48}$.
Then for example, we can take following values for $m$ and an initial value $x_0$.
\begin{center}
    \begin{tabular}{ | l | l | l | l| l|}
    \hline
    $k \pmod {48}  $  & $m$  & $x_0$         \\ \hline
     $8$           &  $19^4$  &$8\cdot13$       \\ \hline
    $12$           &  $20^4$   & $5\cdot17$    \\ \hline
    $20$           &  $2^4$    & $13$          \\ \hline
    $24$           &  $21^4$   & $7\cdot257$   \\ \hline
    $36$           &  $25^4$   & $9\cdot673$   \\ \hline
    $44$           &  $43^4$   & $673$         \\ \hline   
    \end{tabular}
\end{center}
These are easy to check  using a computer.
Note that for these cases, $\gcd(m, G_k)=1$ since a prime divisor of $m$ is either $5$ or congruent to $3 \pmod 4$.
In the above list, we excluded the cases $k\equiv 0, 32 \pmod {48}$.
It seems that there are no small values which satisfy the conditions.
Alternatively, we can further increase the modulus.
Now let us consider it modulo $144$.
Then the remaining cases $k\equiv 0, 32 \pmod {48}$ become $k \equiv 0$, $32$, $48$, $80$, $96$, $128$ $\pmod {144}$.
Then for example, we can take the following values.
\begin{center}
    \begin{tabular}{ | l | l | l | l| l|}
    \hline
    $k \pmod {144}  $  & $m$  & $x_0$               \\ \hline
     $32$           &  $6^4$   &$73$       \\ \hline
    $48$           &  $18^4$   & $2\cdot3\cdot19$        \\ \hline
    $80$           &  $5^4$    & $13$          \\ \hline
    $96$           &  $99^4$   & $3\cdot433$   \\ \hline
    $128$           &  $65^4$   & $2\cdot13$   \\ \hline
        \end{tabular}
\end{center}
Again, we excluded the case when $k\equiv0 \pmod {144}$.
Here again, note that for these cases $\gcd(m, G_k)=1$ since a prime divisor of $m$ is either $5$ or congruent to $3 \pmod 4$.
If we allow a larger modulus, then we might find a  set of initial values for every $k$.
(We want an initial value when $k\equiv0 \pmod {144}$.)

Once we have set an initial value, then the algorithm to check the primality of $H_k$ is the same as the algorithm for $G_k$, simply replace the initial value and replace $G_k$ by $H_k$.

\providecommand{\bysame}{\leavevmode\hbox to3em{\hrulefill}\thinspace}
\providecommand{\MR}{\relax\ifhmode\unskip\space\fi MR }
\providecommand{\MRhref}[2]{%
  \href{http://www.ams.org/mathscinet-getitem?mr=#1}{#2}
}
\providecommand{\href}[2]{#2}

\end{document}